\documentclass[a4paper, twoside, 12pt, final]{amsart}
\usepackage[OT1]{fontenc}
\usepackage[latin1]{inputenc}
\usepackage{amssymb,latexsym}
\usepackage{type1cm,ae}

\theoremstyle{plain}
\newtheorem{theorem}{Theorem}
\newtheorem{corollary}[theorem]{Corollary}
\newtheorem{proposition}[theorem]{Proposition}
\newtheorem{lemma}[theorem]{Lemma}

\theoremstyle{definition}

\theoremstyle{remark}

\newcommand*{\bR}{\ensuremath{\mathbb{R}}}
\newcommand*{\bN}{\ensuremath{\mathbb{N}}}

\newcommand*{\loc}{\mathrm{loc}}

\newcommand*{\Wert}{\mathord{\mbox{|\kern-1.5pt|\kern-1.5pt|}}}

\DeclareMathOperator{\dist}{dist}
\DeclareMathOperator{\osc}{osc}

\DeclareMathOperator{\diam}{diam}

\def\Xint#1{\mathchoice
   {\XXint\displaystyle\textstyle{#1}}%
   {\XXint\textstyle\scriptstyle{#1}}%
   {\XXint\scriptstyle\scriptscriptstyle{#1}}%
   {\XXint\scriptscriptstyle\scriptscriptstyle{#1}}%
   \!\int}
\def\XXint#1#2#3{{\setbox0=\hbox{$#1{#2#3}{\int}$}
     \vcenter{\hbox{$#2#3$}}\kern-.5\wd0}}
\newcommand{\meanint}{\Xint-}

\title[Generalized dimension distortion under Sobolev mappings]
{Generalized Hausdorff dimension distortion in Euclidean spaces under Sobolev mappings}

\author[T. Rajala]{Tapio Rajala}
\author[A. Zapadinskaya]{Aleksandra Zapadinskaya}
\author[T. Z\"urcher]{Thomas Z\"urcher}
\address{Department of Mathematics and Statistics, University of Jyv\"askyl\"a, P.O. Box 35, Fin-40014 University of Jyv\"askyl\"a, Finland}
\email{tapio.m.rajala@jyu.fi}
\email{aleksandra.zapadinskaya@jyu.fi}
\email{thomas.t.zurcher@jyu.fi}


\thanks{The second author was partially supported by the Academy of Finland, grant no.~120972,
and the third author was supported by the Swiss National Science Foundation.}

\begin{document}
\begin{abstract}
We investigate how the integrability of the derivatives of
Orlicz-Sobolev mappings defined on open subsets of $\bR^n$ affect
the sizes of the images of sets of Hausdorff dimension less than
$n$. We measure the sizes of the image sets in terms of generalized
Hausdorff measures.
\end{abstract}
\maketitle

\section{Introduction}
In this paper, we continue the study of Orlicz-Sobolev mappings by
extending planar dimension distortion results to their
$n$-dimensional analogues. We work with mappings that belong to
the class $W^{1,1}_\loc$, meaning that the components of the
mappings have locally integrable distributional derivatives. In
order to obtain interesting distortion results, we make further
assumptions on the mappings.

The distortion is measured on the level of generalized Hausdorff
measure. The \emph{generalized Hausdorff measure} of a set
$A\subset \bR^n$ is defined as
$$
\mathcal{H}^h(A)=\lim_{\delta\to0}\mathcal{H}^h_\delta(A),
$$
where
$$
\mathcal{H}^h_\delta(A)=\inf\Bigl\{\sum\limits_{i=1}^{\infty}
h(\diam U_i)\colon A\subset\bigcup\limits_{i=1}^{\infty} U_i,\,
\diam U_i\leq\delta\Bigr\}
$$
and $h\colon [0,\infty[ \to [0,\infty[$ is a dimension gauge:
$\lim_{t\to0+}h(t)=h(0)=0$ and $h$ is non-decreasing. In
the special case where $h(t)=t^\alpha$ with some $\alpha\geq0$, we have the usual
 \emph{Hausdorff $\alpha$-dimensional measure}, which we simply denote by
$\mathcal{H}^\alpha$. The \emph{Hausdorff
dimension} $\dim_{\mathcal{H}}A$ of a set $A \subset \bR^n$
is the smallest $\alpha_0\geq0$ such that $\mathcal{H}^\alpha(A)=0$ for any
$\alpha>\alpha_0$.

Let $\Omega\subset\bR^n$ be an open set and $f\colon\Omega\to\bR^n$ a continuous
mapping. We obtain distortion estimates of the following form: If the
mapping $f$ is assumed to be in an appropriate Orlicz-Sobolev class, then we have for every $E\subset\Omega$
\[
\dim_{\mathcal H}(E) < n \Longrightarrow \mathcal H^{h_\gamma}(f(E)) = 0,
\]
where $h_\gamma(t) = t^n \log^\gamma(1/t)$. Estimates of this kind
were established in the plane in~\cite{KZZ1,KZZ,Tapio}. The
parameter $\gamma$ in the dimension gauge is chosen according to
the integrability of the differential of the mapping $f$. The
dependence between the parameter $\gamma$ and the threshold for
the integrability of the differential varies with the global
assumptions on the mapping $f$.

We make the strongest integrability assumptions for mappings that
are only assumed to be continuous. We may then relax the
assumptions if we require that in addition to being continuous,
the mappings are also monotone. Recall that a real valued function
$f\colon\Omega\to\bR$ is said to be \emph{monotone} if for every
ball $B\subset\Omega$ we have
$$
\sup_{\partial B}f=\sup_Bf\,\,\,\,\,\,\,\,\text{ and }\,\,\,\,\,\,\,
\inf_{\partial B}f=\inf_Bf.
$$
A mapping $f\colon\Omega\to\bR^n$ is called monotone, if all its
component functions are monotone. Further relaxation may be made
if we assume that the mappings are homeomorphisms.

We list our results for the dimension distortion under
Orlicz-Sobolev mappings under the three different assumptions in the
following theorem.

\begin{theorem}\label{os} Let $\Omega\subset\bR^n$ be an open set and $f:\Omega\to\bR^n$ a continuous
map in $W^{1,1}_\loc(\Omega,\bR^n)$ with $|Df|^n \log^\lambda(e + |Df|)\in
L^1_\loc(\Omega)$ for some\linebreak$\lambda\in\bR$. Then
\begin{equation*}
\dim_{\mathcal H}(E) < n \Longrightarrow \mathcal
H^{h_\gamma}(f(E)) = 0,
\end{equation*}
if one of the following cases occurs:
\begin{itemize}
\item[(\emph{i})] $\lambda > n-1$ and $\gamma < \lambda - n + 1$.
\item[(\emph{ii})]
$f$ is monotone, $\lambda > 0$, and $\gamma \leq  \lambda$.
\item[(\emph{iii})]
$f$ is a homeomorphism, $f^{-1}\in W_\loc^{1,p}(f(\Omega),\bR^n)$
for some \mbox{$p>n-1$}, $\lambda > -1$, and $\gamma \leq
\lambda+1$.
\end{itemize}
\end{theorem}

Theorem \ref{os} will be proved in the remaining sections of
the paper. The claims (\emph{ii}) and (\emph{iii}) are the
analogues of the results known in the planar
case:~~\cite[Theorem~2]{KZZ} and~\cite[Theorem~1.1]{Tapio}. Their
proofs follow the planar ones. However, we have changed the presentation of the proof of claim~(\emph{ii}) in order to emphasize the common key elements of the proofs of claim~(\emph{i}) and claim~(\emph{ii}). Both of these claims follow from a tailored version of a Rado-Reichelderfer condition. To obtain the condition for claim~(\emph{i}), we rely on an auxiliary result from \cite{KKM99}.

We do not know if the estimates in~(\emph{i}) and~(\emph{ii}) are
sharp. However, by~\cite[Proposition~5.1]{HK03}, see Section~2 in
\cite{KZZ1} as well, given any
$\lambda>0$, we may find a homeomorphism
$f\colon[0,1]^n\to[0,1]^n$ such that $|Df|^n\log^s(e+|Df|)$ is
integrable for all $s<\lambda-1$,  mapping a set of
Hausdorff dimension strictly less than $n$ onto a set of positive
generalized Hausdorff measure with the gauge function
$h(t)=t^n\log^\lambda(1/t)$. This example demonstrates the
sharpness of the estimate in~(\emph{iii}).

Theorem~\ref{os} implies dimension distortion results for mappings
of finite distortion. Recall that a continuous mapping $f \in
W^{1,1}_\loc(\Omega, \bR^n)$ is called a \emph{mapping of finite
distortion}, if its Jacobian $J_f$ is locally integrable and there
exists a measurable function $K \colon \Omega \to [1,\infty[$ such
that
\[
 |Df(x)|^n \le K(x)J_f(x)
\]
at almost every point $x \in \Omega$. If we assume $K$ to be
bounded, we obtain the class of \emph{quasiregular mappings}.
However, weaker assumptions already imply distortion estimates.
For example, the assumption that the function $\exp(\lambda K)$ is
locally integrable for some parameter $\lambda>0$. The mappings
for which this is true are called \emph{mappings of
$\lambda$-exponentially integrable distortion}.
See~\cite{Geh-Väi,gehring,BMT} for dimension distortion results
for quasiconformal mappings and~\cite{HK03,KZZ} for 
generalized dimension distortion estimates for mappings of
exponentially integrable distortion in the plane.

We follow the approach taken in~\cite{KZZ} and prove dimension
distortion results for mappings of $\lambda$-exponentially
integrable distortion using the higher regularity of the weak
derivatives of the mappings. We obtain the higher regularity from \cite{FKZ05}.  In the plane, the sharp regularity is known, see \cite{eero}. A
straightforward  combination of Theorem~\ref{os}
and~\cite[Theorem~1.1]{FKZ05} gives the following result.

\begin{corollary}\label{fd} Let $\Omega\subset\bR^n$ be a domain and $\lambda>0$. There exist
positive constants $c_1$ and $c_2$ depending only on $n$ such that if $f:\Omega\to\bR^n$ is of
$\lambda$-exponentially integrable distortion and satisfies
\begin{itemize}
\item[(\emph{i})]
$\lambda > 1/c_1$ and $\gamma \leq  c_1\lambda - 1$ or
\item[{(\emph{ii})}]
$f$ is homeomorphic and $\gamma \leq c_2\lambda$,
\end{itemize}
the following implication is true:
\begin{equation*}
\dim_{\mathcal H}(E) < n \Longrightarrow \mathcal
H^{h_\gamma}(f(E)) = 0.
\end{equation*}
\end{corollary}

\begin{proof}
By~\cite[Theorem~1.1]{FKZ05} there exists a constant $c_1$ depending only on $n$ such that
\begin{equation}\label{eq:corollary}
|Df|^n\log^{c_1\lambda-1}(e+|Df|)\in L^1_\loc(\Omega)
\end{equation}
for each mapping $f$ satisfying the assumptions of the corollary.

By a result in~\cite{vodopis} saying that mappings of finite
distortion of the class $W^{1,n}_\loc(\Omega,\bR^n)$ are continuous and monotone,
claim (\emph{i}) follows from Theorem~\ref{os} (\emph{ii}).

For the proof of (\emph{ii}), we note that \eqref{eq:corollary}
implies that $|Df|$ is in the Lorentz space $L^{n-1,1}_\textnormal{loc}(\Omega)$, see
\cite[V.3]{Lorentz}. The assumptions of Theorem~4.1
in~\cite{HeKoMa} are fulfilled in our settings, giving us
$f^{-1}\in W^{1,n}_\loc(f(\Omega),\Omega)$. This together
with~\eqref{eq:corollary} allows us to apply (\emph{iii}) of
Theorem~\ref{os}, which concludes the proof.
\end{proof}
The estimate in Corollary~\ref{fd} (\emph{ii}) is sharp modulo the
constant. Indeed, using again~\cite[Proposition~5.1]{HK03}, for
any given $0<\varepsilon<\lambda$, we find a homeomorphism
$f\colon[0,1]^n\to[0,1]^n$, having
$\nobreak{(\lambda-\varepsilon)}$\nobreakdash-exponentially
integrable distortion and mapping a set of Hausdorff dimension
strictly less than $n$ onto a set of positive generalized
Hausdorff measure with the gauge function
$h(t)=t^n\log^\lambda(1/t)$. This example makes us believe that
(\emph{ii}) in Corollary~\ref{fd} holds with $c_2=1$. However, we
do not know if the estimate in (\emph{i}) is sharp.

In the next section, we prove a result that we will use for proving the cases (\emph{i}) and (\emph{ii}) in Theorem~\ref{os}. Afterwards, three more sections follow --- one for each case of Theorem~\ref{os}.

\section{Rado-Reichelderfer condition}
A version of the following inequality~\eqref{RReq} was first mentioned by Rado and Reichelderfer, \cite{RadoReichelderfer}.

\begin{proposition}\label{RR}
Suppose that $\Omega\subset
\bR^n$ is an open set, $f\colon \Omega\to\bR^n$ a continuous mapping, and $\gamma>0$. Assume further that there is a function $\rho \colon \Omega \to [0,\infty]$ and constants $C\ge1$ and $r_0>0$ such that
\begin{equation*}
\rho\log^\gamma(e+\rho)\in L^1_\textnormal{loc}(\Omega)
\end{equation*}
and
\begin{equation}\label{RReq}
(\diam f(B(x,r)))^n\leq \int_{B(x,r)} \rho
\end{equation}
for all $0<r<r_0$ and $x\in \Omega$ taken so that $B(x,Cr) \subset \Omega$. Then
\begin{equation*}
\dim_{\mathcal{H}}(E)<n\Longrightarrow
\mathcal{H}^{h_\gamma}(f(E))=0.
\end{equation*}
\end{proposition}

\begin{proof}
Let $E\subset \Omega$ be such that $\dim_{\mathcal{H}}(E)<n$. By the $\sigma$\nobreakdash-additivity of the Hausdorff measure, we may assume without
loss of generality that $E \subset\subset \Omega$. Fix $0<\varepsilon<\min\{1,e^{-\lambda n}\}$. By the absolute continuity of the integral, we may find
\mbox{$\delta\in]0,\min\{\frac{1}{4},e^{-\frac{2\lambda}{n}},r_0\}[$} such
that
\begin{equation}\label{intbound}
\int_{A}\rho\log^\gamma(e+\rho)<\varepsilon<\min\{1,e^{-\lambda n}\},
\end{equation}
whenever $A\subset\Omega$ is such that $\mathcal L^n (A)<\delta$, which is for example the case if $A$ is a ball of radius less than $\delta$. Note that 
\[
\rho(x) \le \rho(x)\log^\gamma(e + \rho(x)), 
\]
so we obtain from \eqref{RReq}
\begin{equation}\label{smallballs}
\diam f(B(x,r)) \le \left(\int_{B(x,r)}\rho \right)^{1/n}<\varepsilon^{1/n} < e^{-\lambda}
\end{equation}
for $B(x,r)\subset\Omega$ with $0<r<\delta$.

Let us describe how to choose a suitable cover for $E$ that gives rise
to an eligible cover of $f(E)$. First, we take $\alpha\in]\max\{\frac{n}{2},\dim_{\mathcal{H}}E\},n[$. Notice that
there exists $t_0>0$ such that
\begin{equation}\label{logestimate}
 t^n \le t^\alpha \log^\gamma\left(\frac{1}{t}\right) \le t^{(\alpha +\dim_{\mathcal{H}}E)/2}
\end{equation}
for all $0<t<t_0$. Since $\mathcal{H}^{(\alpha
+\dim_{\mathcal{H}}E)/2}(E) = 0$, we may find for
every\linebreak$\varepsilon'>0$ a covering of $E$ with balls
$\{B_i'\}_{i=1}^\infty$ of diameter less
than\linebreak$\min\{t_0,\delta, \frac{1}{C}\dist(E, \mathbb{R}^n \setminus \Omega)\}$ so that
\[
 \sum_{i=1}^\infty(\diam B_i')^{(\alpha +\dim_{\mathcal{H}}E)/2} < \varepsilon'.
\]
From the cover $\{B'_i\}$, we move to a more suitable cover. This is done by defining a new collection of balls
\[
 \mathcal{B}' = \big\{B(x,r) : x \in B_i'\cap E \text{ and }2r = \diam B_i'\text{ for some }i \in \bN\big\}
\]
and applying the Besicovitch covering theorem to it. This gives a
constant $0<N<\infty$, depending only on $n$, and a covering of the set
$E$ with balls $\mathcal{B}=\{B_j\}_{j=1}^\infty \subset \mathcal{B}'$
of diameters less than $\delta$ so that $\sum_j \chi_{B_j}(x) \le N$ for every $x \in \Omega$.

Now, for any ball $B'_i$ from the original cover, all the balls $B(x,r)$ in $\mathcal{B}$ with $x \in B'_i$ and $2r = \diam B_i'$ contain the center of the ball $B'_i$, hence there are at most $N$ such balls. Therefore
\[
 \sum_{j=1}^\infty(\diam B_j)^{(\alpha +\dim_{\mathcal{H}}E)/2} \le N\sum_{i=1}^\infty(\diam B_i')^{(\alpha +\dim_{\mathcal{H}}E)/2} < N\varepsilon'.
\]
Now, by taking $\varepsilon'$ small enough, we have
 $\mathcal L^n(\bigcup_j B_j)<\delta$ and by \eqref{logestimate}
$$
\sum_{j=1}^\infty (\diam B_j)^\alpha\log^\gamma\Bigl(\frac{1}{\diam
B_j}\Bigr)<\varepsilon.
$$

Let us show that $\mathcal{H}_\delta^{h_\gamma}(f(E))=0$. We use the monotonicity of $t\log^\lambda(1/t)$ for $t\in
]0,e^{-\lambda}[$ and estimate~\eqref{smallballs} together with the Rado-Reichelderfer condition~\eqref{RReq} to obtain
\begin{align}\label{StartEq}
(\diam f(B))^n\log^\gamma\Bigl(\frac{1}{(\diam f(B))^n}\Bigr)\leq
\int_B \rho\cdot \log^\gamma\Bigl(\frac{1}{\int_B \rho}\Bigr).
\end{align}
We consider two cases. Let
\begin{align*}
\mathcal{B}_1:=\Big\{B_j\in \mathcal{B}:\, \int_{B_j}\rho\leq (\diam B_j)^\alpha\Big\},\\
\mathcal{B}_2:=\Big\{B_j\in \mathcal{B}:\, \int_{B_j}\rho > (\diam B_j)^\alpha\Big\}.
\end{align*}
For $B\in \mathcal{B}_1$, we use again the monotonicity of
$t\log^\lambda(1/t)$ to obtain from inequality \eqref{StartEq}
\begin{align*}
(\diam f(B))^n\log^\gamma\Bigl(\frac{1}{(\diam f(B))^n}\Bigr)\leq
\alpha^\gamma(\diam B)^\alpha \log^\gamma\Bigl(\frac{1}{\diam B}\Bigr).
\end{align*}
In case that $B \in \mathcal{B}_2$, we conclude from \eqref{StartEq} that
\begin{align*}
(\diam f(B))^n\log^\gamma\Bigl(\frac{1}{(\diam f(B))^n}\Bigr)
\leq 
\alpha^\gamma\int_B \rho\cdot \log^\gamma\Bigl(\frac{1}{\diam
B}\Bigr).
\end{align*}
We split the integral in two parts. In order to continue, we  set
\[
A_B:=\Big\{y\in B:\, \rho(y)\le\frac{1}{(\diam B)^{n-\alpha}}\Big\}
\]
and obtain
\begin{align*}
\int_{A_B} \rho(y)&\cdot  \log^\gamma\Bigl(\frac{1}{\diam B}\Bigr)\, dy\\
&\leq
\frac{\omega_n}{(\diam
B)^{n-\alpha}}(\diam B)^n\log^\gamma\Bigl(\frac{1}{\diam B}\Bigr)\\
&= \omega_n(\diam B)^\alpha\log^\gamma\Bigl(\frac{1}{\diam B}\Bigr),
\end{align*}
where $\omega_n$ denotes the $n$\nobreakdash-dimensional Lebesgue measure  of an $n$-dimensional unit ball. Now, we focus on $B\setminus A_B$:
\begin{align*}
\int_{B\setminus A_B}&\rho(y)\log^\gamma\Bigl(\frac{1}{(\diam
B)^\alpha}\Bigr)\, dy\\
=&\int_{B\setminus A_B} \rho(y)\Bigl(\frac{\alpha}{n-\alpha}\Bigr)^\gamma
\log^\gamma\Bigl(\frac{1}{(\diam B)^{n-\alpha}}\Bigr)\, dy \\
\leq& \Bigl(\frac{\alpha}{n-\alpha}\Bigr)^\gamma\int_{B\setminus A_B}
\rho(y)
\log^\gamma(e+\rho(y))\, dy.
\end{align*}
Hence, keeping \eqref{intbound} in mind, we obtain the following upper bound:
\begin{align*}
\sum\limits_{B\in\mathcal
B_2}&\int_{B\setminus A_B}\rho(y)\log^\gamma\Bigl(\frac{1}{(\diam
B)^\alpha}\Bigr)\, dy\\
&\leq\Bigl(\frac{\alpha}{n-\alpha}\Bigr)^\gamma N\int_{\bigcup_i B_i}\rho(y)
\log^\gamma(e+\rho(y))\, dy\\
&<N\Bigl(\frac{\alpha}{n-\alpha}\Bigr)^\gamma
\varepsilon.
\end{align*}
Thus, we estimate
\begin{align*}
 \mathcal H^h_{\varepsilon^{1/n}}(f(E))\leq&n^{-\gamma}\sum_{j=1}^\infty (\diam
f(B_j))^n\log^\gamma\Bigl(\frac{1}{(\diam
f(B_j))^n}\Bigr)\\
\leq&\omega_n
n^{-\gamma}\alpha^\gamma\Bigl(1+N\Bigl(\frac{\alpha}{n-\alpha}
\Bigr)^\gamma \Bigr)\varepsilon.
\end{align*}
Letting $\varepsilon$ go to zero concludes the proof.
\end{proof}

\section{Proof of the continuous case}
Suppose that $\varphi$ is a positive function on the interval
$]0,\infty[$. We write
\begin{equation*}
F_\varphi(s)=
\begin{cases}
s\varphi^{\frac{1}{n}-1}(s) & s>0,\\
0 & s=0.
\end{cases}
\end{equation*}

Our proof of Theorem \ref{os} (i) is based on~\cite[Theorem~3.2]{KKM99}, which, in the case of
continuous mappings, states the following.
\begin{theorem}\label{Lorentz}
Let $u\in W^{1,1}_\loc(\Omega,\bR)$ be a continuous
function. Let further $\varphi$ be a positive, non-increasing
function on $]0,\infty[$. Suppose that
\begin{equation*}
\int_\Omega F_\varphi(|Du(x)|)\, dx<\infty
\end{equation*}
and
\begin{equation*}
\int_0^\infty \varphi^{1/n}(t)\, dt<\infty.
\end{equation*}
Then
\begin{equation*}
(\osc_{B(x,r)} u)^n = \int_{B(x,r)}
\frac{2^{n(n+2)}}{n\omega_n}\biggl(\int_0^\infty
\varphi^{1/n}(t)\, dt\biggr)^{n-1} F_\varphi(|Du(x)|)\, dx.
\end{equation*}
\end{theorem}

\begin{proof}[Proof of Theorem~\ref{os}~(i)]
Let us fix $\gamma\in]0,\lambda+1-n[$ and prove that $\mathcal
H^h(f(E))=0$ with $h(t)=t^n\log^\gamma(1/t)$. We want to verify the Rado-Reichelderfer condition~\eqref{RReq} by applying Theorem~\ref{Lorentz}. Hence, we need to choose a function $\varphi\colon
]0,\infty[\to \bR$ so that $F_\varphi$ matches the integrability
condition specified in Theorem~\ref{os}. We choose \mbox{
$\tilde{\lambda}\in]n-1, \lambda-\gamma[$} and set
\begin{equation*}
\varphi(s)=
\begin{cases}
\log^\frac{\tilde{\lambda} n}{1-n}(e+1) & 0<s<1,\\
s^{-n}\log^\frac{\tilde{\lambda} n}{1-n}(e+s) & s\geq 1.
\end{cases}
\end{equation*}
This is a suitable choice for $\varphi$ since then $F_\varphi$ is defined as
\begin{equation*}
F_\varphi(s)=
\begin{cases}
0 & s=0,\\
s\log^{\tilde{\lambda}}(e+1) & 0<s\leq 1,\\
s^n\log^{\tilde{\lambda}}(e+s) & 1<s
\end{cases}
\end{equation*}
giving
\begin{equation*}
\int_\Omega F_\varphi(|Df(x)|)\, dx <\infty
\end{equation*}
by assumption. Clearly, $\varphi$ satisfies the requirements in
Theorem~\ref{Lorentz}.

We conclude from Theorem \ref{Lorentz} that there exists a constant $C>0$ depending only on $n$ and $\tilde{\lambda}$ so that for every $x \in \Omega$ and $r > 0$ with $B(x,r) \subset \Omega$, the estimate 
\begin{align}\label{mainineq}
(\diam f(B(x,r)))^n\leq& \Bigl(\sum\limits_{i=1}^n\osc^2_{B(x,r)}f_i\Bigr)^{n/2}\\
\leq&
\Bigl(\sum\limits_{i=1}^n\Bigl(C\int_{B(x,r)}F_\varphi(|Df_i|)\Bigr)^{2/n}\Bigr)^{n/2
}\notag\\
\leq& \,Cn^{n/2}\int_{B(x,r)}F_\varphi(|Df|)\notag
\end{align}
holds. We define $\rho$ as $Cn^{n/2}F_\varphi(|Df|)$. 
We further notice that 
\[
\rho\log^\gamma(e+\rho)
\]
is locally  integrable in $\Omega$
due to the choice of parameters and the given integrability of
$|Df|$; indeed, \mbox{$s\log^\gamma(e+s)$} behaves
asymptotically like $s\log^{\tilde{\lambda}+\gamma}(e+s)$ when
$s$ is large, and we have \mbox{$\tilde{\lambda}+\gamma<\lambda$}. Having verified the Rado-Reichelderfer condition as specified in Proposition~\ref{RR}, we conclude the proof.
\end{proof}

\section{Proof of the continuous and monotone case}
We use the symbol $A(x,r,R)$ to denote the closed annulus with center at $x$
and radii $r$ and $R$:
$$
A(x,r,R)=\{y\in\bR^n\colon r\leq|x-y|\leq R\}.
$$
The proof we present here repeats the main steps of its planar
analogue ~\cite[Theorem~2]{KZZ}. We start by proving an estimate
on the sizes of the images of balls.
\begin{lemma}\label{lem:uboj}
Let $\Omega$ be a domain in $\bR^n$ and $f\colon\Omega\to
f(\Omega)$ be a continuous and monotone mapping of the Sobolev class
$W^{1,p}_\loc(\Omega, \bR^n)$, where $n-1<p<n$. Then there exists a
constant $C$ depending only on $n$ and $p$ such that for every ball
$B(x,2r)\subset\subset\Omega$, the following estimate holds
\begin{equation}\label{estimateAnnulus}
 \diam f(B(x,r))\leq Cr^{1-\frac{n}{p}}
\Bigl(\int_{A(x,r,2r)}|Df|^p\Bigr)^{1/p}.
\end{equation}
\end{lemma}
\begin{proof} Fix a ball $B(x,r)$. Without loss of generality, we may
  assume that $B(x,2r)\subset\subset\Omega'$ for some domain $\Omega'$
  in $\Omega$ and that $f\in W^{1,p}(\Omega',\bR^n)$.  Using the monotonicity of $f$, we obtain
\begin{align*}
\diam f(B(x,r))&\leq(\sum_{i=1}^{n}(\osc_{B(x,r)} f_i)^2)^{1/2}\\
&\leq\sum\limits_{i=1}^{n}\osc\limits_{B(x,r)} f_i\leq\sum\limits_{i=1}^{n}\osc\limits_{S(x,t)} f_i 
\end{align*}
for all $t\in]r,2r[$. Since $f_i$ is in $W^{1,p}(\Omega',
\bR)$, $i=1,\ldots,n$, Theorem~5.16 on p.~121 in \cite{Degree} gives us $f_i\in W^{1,p}(S(x,t))$, $i=1,\ldots,n$, for
almost every $t\in]r,2r[$. See the same page for the
definition of the class $W^{1,p}(X)$ for a $C^\infty$ paracompact
manifold $X$ such as a sphere and for Theorem~5.15, which
gives us
$$
\osc\limits_{S(x,t)} f_i\leq
Ct^{1-\frac{n-1}{p}}\Bigl(\int_{S(x,t)}|Df_i|^p\Bigr)^{1/p}
$$
for almost every $t\in [r,2r]$, $i=1,\ldots,n$. Hence
$$
\diam f(B(x,r))\leq Ct^{1-\frac{n-1}{p}}\sum\limits_{i=1}^{n}\Bigl(\int_{S(x,t)}|Df_i|^p\Bigr)^{1/p}
$$
for almost every $t\in[r,2r]$. Integrating the last inequality over $[r,2r]$ with respect to $t$
and using H\"older's inequality, we obtain
\begin{align*}
 r\diam &f(B(x,r))\leq C\int_{[r,2r]}t^{1-\frac{n-1}{p}}
\sum\limits_{i=1}^{n}\Bigl(\int_{S(x,t)}|Df_i|^p\Bigr)^{1/p}dt\\
&\leq C\sum\limits_{i=1}^{n}\Bigl(\int_{[r,2r]}\int_{S(x,t)}|Df_i|^p\,dt\Bigr)^{1/p}
\Bigl(\int_{[r,2r]}t^{\frac{p-n+1}{p-1}}\Bigr)^{\frac{p-1}{p}}\\
&\leq Cr^{2-\frac{n}{p}}\sum\limits_{i=1}^{n}\Bigl(\int_{A(x,r,2r)}|Df_i|^p\Bigr)^{1/p}
\\&\leq Cnr^{2-\frac{n}{p}}\Bigl(\int_{A(x,r,2r)}|Df|^p\Bigr)^{1/p}.
\end{align*}
\end{proof}
In comparison with the Rado-Reichelderfer condition~\eqref{RReq}, we
integrate in \eqref{estimateAnnulus} over the annulus $A(x,r,2r)$
instead of the  ball $B(x,r)$. Our aim is to bound the integral over the annulus by an integral over the ball. In order to do so, we have to change the integrand. More precisely, we replace it by a maximal operator.

Assume that
$\Omega\subset\bR^n$ is a cube and $h\colon \Omega \to \bR$ is a
non-negative and integrable function. The maximal operator
$\mathcal{M}_\Omega$ at a point $x \in \Omega$ is defined by
\begin{equation*}
\mathcal{M}_\Omega h(x)=\sup\Bigl\{\meanint_{Q} h\, dy\colon\, x\in Q\subset
\Omega\Bigr\},
\end{equation*}
where the supremum is taken over all subcubes of $\Omega$ containing the point $x$.

In order to continue, we use the following lemma. It was proved in~\cite[Lemma
2]{KZZ} in the planar case. However, the same proof works in higher dimensions.
\begin{lemma}\label{lemma2}
 Let $\Omega\subset \bR^n$ be a cube, $f\colon\Omega\to \bR^n$ be a mapping of the Sobolev class
$W^{1,p}_\loc(\Omega, \bR^n)$, where $n-1<p<n$. Then there exists a constant $C$ depending only on $n$
such that the inequality
$$
\int_{A(x,r,2r)}|Df|^p\leq
C\int_{B\left(x,r\right)}\mathcal{M}_\Omega|Df|^p
$$
holds for all $x\in\Omega$ and $r>0$, such that $B(x,4\sqrt{n}r)\subset\subset\Omega$.
\end{lemma}
We will also use the following auxiliary result, which was proved in \cite[Lemma 5.1]{Limits}.
\begin{lemma}\label{lemmafromLimits}
 Suppose $A \colon [0,\infty[ \to [0,\infty[$ is increasing and $\Phi(t) = A(t)t^p$ for some $p > 1$. Then there exists a constant $C>0$ depending only on $n$ and $p$ such that
\[
 \int_\Omega \Phi(M_\Omega h) \le C \int_\Omega \Phi(Ch).
\]

\end{lemma}

Now, we are ready to combine the preceding 3 lemmas.

\begin{proposition}\label{prop:1}
Let $\Omega\subset \bR^n$ be a cube and $f\colon \Omega\to
f(\Omega)$ be a continuous and monotone mapping in $W^{1,1}\left(\Omega,\bR^n\right)$ satisfying
$$|{Df}|^n\log^\lambda\left(e+|{Df}|\right)\in
L^1\left(\Omega\right)$$
for some $\lambda>0$. Then there exists a function $\rho\in L^1(\Omega)$ such that
\begin{equation}\label{intOfRho}
\rho\log^\lambda(e+\rho)\in L^1(\Omega)
\end{equation}
and
\begin{equation}\label{RRConclusion}
(\diam f(B\left(x,r\right)))^n\leq \int_{B\left(x,r\right)}\rho,
\end{equation}
for all $x\in \Omega$ and $r>0$ such that $B(x,4\sqrt{n}r)\subset\subset\Omega$.
\end{proposition}
\begin{proof}
The proof is similar to the one in~\cite{KZZ}. Fix $x\in \Omega$ and $r>0$ such that $B(x,4\sqrt{n}r)\subset\subset\Omega$. Fix also $p\in]n-1,n[$. In Lemma~\ref{lem:uboj}, we deduced the inequality
$$
 \diam f(B(x,r))\leq Cr^{1-\frac{n}{p}}
\Bigl(\int_{A(x,r,2r)}|Df|^p\Bigr)^{1/p}
$$
with a constant $C$ depending only on $n$ and $p$. In what follows, the constant $C$, still depending only on $n$ and $p$, may change its value
from occurrence to occurrence.
The combination of the inequality above with the upper bound found in Lemma~\ref{lemma2} implies
$$
\diam f(B(x,r))\leq Cr^{1-\frac{n}{p}}
\Bigl(\int_{B(x,r)}\mathcal{M}_\Omega(|Df|^p)\Bigr)^{1/p}.
$$
We continue the estimation with the
H\"older inequality:
\begin{align}\label{eq:keyineq}
(\diam f(B(x,r)))^n&\leq Cr^{n-\frac{n^2}{p}}
\Bigl(\int_{B(x,r)}\mathcal{M}_\Omega(|Df|^p)\Bigr)^{n/p}\\
&\leq C\int_{B(x,r)}\mathcal{M}^{n/p}_\Omega(|Df|^p).\notag
\end{align}
In order to finish, we have to verify that
\begin{equation*} 
\rho(x):=\mathcal M^{n/p}_\Omega(|Df|^p)(x)
\end{equation*}
fulfills~\eqref{intOfRho}. As in Lemma~3 in \cite{KZZ},  we prove it as an
application of Lemma~\ref{lemmafromLimits}. We let $h=|Df|^p$ and $\Phi(t)=t^{n/p}\log^\lambda(e+t^{n/p})$ and use the fact that $n/p>1$. We get
\begin{align*}
\int_\Omega&\mathcal{M}_\Omega^{n/p}(|Df|^p)
\log^\lambda\bigl(e+\mathcal{M}_\Omega^{n/p}(|Df|^p)\bigr)\\
&\leq C\int_\Omega|Df|^n
\log^\lambda(e+C|Df|^n)\\
&\leq C\int_\Omega|Df|^n
\log^\lambda(e+|Df|)<\infty.
\end{align*}
\end{proof}
\begin{proof}[Proof of Theorem~\ref{os}~(ii)]
By the $\sigma$\nobreakdash-additivity of the Hausdorff measure, we may assume that $\Omega$ is a cube. In Proposition~\ref{prop:1}, we have verified the Rado-Reichelderfer condition~\eqref{RRConclusion}. In Proposition~\ref{RR}, we have shown that this condition is sufficient under the given integrability of $|Df|$ to conclude the proof.  
\end{proof}

\section{Proof of the homeomorphic case}

For a set $V\subset\bR^n$ and a number $\delta>0$, $V+\delta$ denotes the set
$\{y\in\bR^n\,:\dist(y,V)<\delta\}$.

Without loss of generality, we assume for the rest of this section
that $\Omega$, in addition to being open, is connected. The first
step is an analogue of~\cite[Lemma~3.1]{KZZ1}.
\begin{lemma}\label{lemma:inverse}
 Let $f\colon\Omega\to f(\Omega)\subset\bR^n$ be a homeomorphism such that $f^{-1}\in
W^{1,p}_\loc(\Omega,\bR^n)$ for some $p\in]n-1,n[$. Then there exists a set $F\subset
f(\Omega)$ such that $\mathcal{H}^{n-\frac{p}{2}}(F)=0$ and for all $y\in
f(\Omega)\setminus F$ there exist constants $C_y>0$ and $r_y>0$
such that
\begin{equation}\label{ineq:desired}
\diam(f^{-1}(B(y,r)))\leq C_yr^{1/2},
\end{equation}
for all $0<r<r_y$.
\end{lemma}
\begin{proof}
 Let us apply Lemma~\ref{lem:uboj} to the mapping $f^{-1}$. We obtain
$$
\diam f^{-1}(B(y,r))\leq Cr^{1-\frac{n}{p}}
\sum\limits_{i=1}^{n}\Bigl(\int_{B(y,2r)}|Df^{-1}_i|^p\Bigr)^{1/p}
$$
for all $y\in f(\Omega)$ and $r>0$ such that
$B(y,2r)\subset\subset f(\Omega)$, where the constant $C$ depends only on $n$ and $p$. So, the desired inequality
holds for all $y\in f(\Omega)$ for which
\begin{equation}\label{eq:bad}
r^{\frac{p}{2}-n}\int_{B(y,2r)}|Df^{-1}_i|^p<M_y
\end{equation}
is valid for all $i=1,\ldots,n$, for all small enough $r>0$ and for some constant $M_y$, depending on $y$. Let $F_1$ be the set of those
$y$ for which \eqref{eq:bad} does not hold with $i=1$. Let
$K\subset f(\Omega)$ be a compact set and fix some $\delta>0$ such
that $\dist(K,\partial f(\Omega))>\delta$. For every
$k\in\mathbb{N}$ and every $y$ in $F_1\cap K$, there exists
$r_{k,y}<\delta/10$ such that
$\int_{B(y,2r_{k,y})}|Df^{-1}_1|^p\geq
k(r_{k,y})^{n-\frac{p}{2}}$. Consider the collection of balls
$\mathcal{B}_k=\{B(y,2r_{k,y})\colon y\in F_1\cap K\}$ for every
$k\in\bN $. Using Vitali's covering theorem, we obtain for every
$k\in\mathbb{N}$ a countable subcollection of disjoint balls
$B_{k,j}$, $j=1,2,\ldots$, centered in $F_1\cap K$, having radii
$2r_{k,j}<\delta/5$, and with $\bigcup\limits_{j=1}^\infty5B_{k,j}$
covering $F_1\cap K$. We have
\begin{align*}
\mathcal{H}^{n-\frac{p}{2}}_\delta(F_1\cap
K)&\leq\sum\limits_{j=1}^{\infty}(20r_{k,j})^{n-\frac{p}{2}}\leq\frac{20^{n-\frac{p}{2}}}{k}
\sum\limits_{j=1}^{\infty}\int_{B_{k,j}}|Df^{-1}_1|^p\\&
\leq\frac{20^{n-\frac{p}{2}}}{k}\int_{K+\delta/5}|Df^{-1}_1|^p
\end{align*}
for all $k\in\mathbb N$. Letting $k\to\infty$ and
$\delta\to0$, we obtain \mbox{$\mathcal{H}^{n-\frac{p}{2}}(F_1\cap K)=0$}.
\end{proof}
\begin{proof}[Proof of Theorem~\ref{os}~(iii)]
Since $f$ is a homeomorphism in $W^{1,p}_\loc(\Omega,\bR^n)$ for all $0<p<n$, we know by~\cite{STH} that the Jacobian of $f$ is
either non-positive or non-negative almost everywhere
(we assumed $\Omega$ to be connected). Thus, we may assume that
$J_f\geq 0$ a.e. in $\Omega$. The rest of the proof goes as in~\cite{Tapio}
(Lemma~2.2 and the proof of Theorem~1.1). 

Let us outline the proof. First, we cover the set $E$ with a countable collection of sets $E_j$ so that in each set the constants $C_y$ and $r_y$ of Lemma~\ref{lemma:inverse} are fixed. On each of the covering sets, we use the fact that the dimension of the set is less than $n$ to obtain a nice covering by small balls. Next, we estimate the size of the image of the sets $E_j$ by the integral of the Jacobian of $f$. For this purpose, we enlarge the covering balls, use the $5r$\nobreakdash-covering theorem and Lemma~\ref{lemma:inverse}. The claim follows from the higher integrability of the Jacobian of $f$, which is obtained from the assumption on the integrability of $|Df|$ by \cite[Corollary~9.1]{ivaver}.
\end{proof}

\subsection*{Acknowledgments} The authors thank Pekka Koskela for suggesting this problem.

\bibliographystyle{amsplain}
\bibliography{jyv7}
\end{document}